\documentclass[12pt,a4paper]{article}
\usepackage{mathrsfs}
\usepackage{amssymb}
\usepackage{amsmath}
\usepackage{amsfonts}
\usepackage{amstext}
\usepackage{amsthm}
\usepackage{graphicx}
\usepackage[top=1in, bottom=1in,left=1in,right=1in]{geometry}

\def \qed {\hfill \vrule height6pt width 6pt depth 0pt}

\begin{document}

\title{Unilateral global bifurcation and nodal solutions for the $p$-Laplacian with sign-changing weight
\thanks{Research supported by the NSFC (No.11061030, No.10971087).}}
\author{{\small  Guowei Dai\thanks{Corresponding author. Tel: +86 931
7971297.\newline
\text{\quad\,\,\, E-mail address}: daiguowei@nwnu.edu.cn (G. Dai).%
}\, \, Ruyun Ma
} \\
%EndAName
{\small Department of Mathematics, Northwest Normal
University, Lanzhou, 730070, PR China}\\
}
\date{}
\maketitle

\begin{abstract}
In this paper, we shall establish a Dancer-type unilateral global bifurcation result for a class of quasilinear elliptic problems with sign-changing weight.
Under some natural hypotheses on perturbation function,
we show that $\left(\mu_k^\nu(p),0\right)$ is a bifurcation
point of the above problems and there are two distinct unbounded continua,
$\left(\mathcal{C}_{k}^\nu\right)^+$ and $\left(\mathcal{C}_{k}^\nu\right)^-$,
consisting of the bifurcation branch $\mathcal{C}_{k}^\nu$ from $\left(\mu_k^\nu(p), 0\right)$,
where $\mu_k^\nu(p)$ is the $k$-th positive or negative eigenvalue of the linear problem corresponding to the above problems,
$\nu\in\{+,-\}$. As the applications of the above unilateral global bifurcation result, we study the existence of nodal solutions
for a class of quasilinear elliptic problems with sign-changing weight. Moreover, based on the bifurcation result of Dr\'{a}bek and Huang (1997) [\ref{DH}], we study the existence of one-sign solutions for a class of high dimensional quasilinear elliptic problems with sign-changing weight.
\end{abstract}\textbf{\ }

\numberwithin{equation}{section}

\numberwithin{equation}{section}

\section{\small Introduction}

\quad\, In [\ref{R2}], Rabinowitz established Rabinowitz's unilateral
global bifurcation theory. However,
as pointed out by Dancer [\ref{D1}, \ref{D2}] and L\'{o}pez-G\'{o}mez [\ref{L1}],
the proofs of these theorems contain gaps. Fortunately, Dancer [\ref{D1}] gave a corrected version of unilateral global bifurcation
theorem. In 1997, Dr\'{a}bek and Huang
[\ref{DH}] proved a Dancer-type bifurcation theorem (Theorem 4.5, [\ref{DH}]) in which the continua
bifurcated from the principle eigenvalue for a high dimension $p$-Laplacian problem with sign-changing weight in $\mathbb{R}^N$.
However, no any information on the high eigenvalue for $p$-Laplacian problem with sign-changing weight. For the case of definite weight,
Dai and Ma [\ref{DM}] established a Dancer-type unilateral global bifurcation result for one-dimensional $p$-Laplacian.
In [\ref{GT}], Girg and Tak\'{a}\v{c} proved a Dancer-type bifurcation theorem for a high dimensional $p$-Laplacian equation.

It is the main purpose of this paper to establish a similar result to Dai and Ma's
about the continua of radial solutions for the following $N$-dimensional $p$-Laplacian
problem on the unit ball of $\mathbb{R}^N$ with $N\geq1$ and $1<p<+\infty$,
\begin{equation}\label{egB1}
\left\{
\begin{array}{l}
-\text{div}\left(\varphi_p(\nabla u)\right)=\mu m(x)\varphi_p(u)+g(x,u;\mu), \ \ \  \ \text{in}\,\,
B,\\
u=0,\quad\quad\quad\quad\quad\quad\quad\quad\quad\quad\quad\quad\quad\quad\quad\quad\quad \text{on}\,\,\partial B,
\end{array}
\right.
\end{equation}
where $B$ is the unit open ball of $\mathbb{R}^N$, $\varphi_p(s)=\vert s\vert^{p-2}s$, $m\in M(B)$ is a sign-changing
function with
\begin{equation}
M(B)=\left\{m\in C(\overline{B})\,\,\text{is radially symmetric\,}\big|\text{meas}\{x\in B,m(x)>0\}\neq0\right\},\nonumber
\end{equation}
$g:B\times \mathbb{R}\times\mathbb{R}\rightarrow\mathbb{R}$ satisfies the Carath\'{e}odory condition in the first
two variables and radially symmetric with respect to $x$.

It is clear that the radial solutions of (\ref{egB1}) is
equivalent to the solutions of the following problem
\begin{equation}\label{eg1}
\left\{
\begin{array}{l}
-\left(r^{N-1}\varphi_p(u')\right)'=\mu m(r)r^{N-1}\varphi_p(u)+r^{N-1}g(r,u;\mu), \ \ \text{a.e.}\  \ r\in
(0,1),\\
u'(0)=u(1)=0,
\end{array}
\right.
\end{equation}
where $r=\vert x\vert$ with $x\in B$, $m\in M(I)$ is a sign-changing with $I=(0,1)$ and
\begin{equation}
M(I)=\left\{m\in C(\overline{I})\,\,\text{is radially symmetric\,}\big|\text{meas}\{r\in I,m(r)>0\}\neq0\right\}.\nonumber
\end{equation}
We also assume the perturbation function $g$ satisfies the following hypotheses:
\begin{equation}\label{hg1}
\lim_{s\rightarrow0}\frac{g(r,s;\mu)}{\vert s\vert^{p-1}}=0
\end{equation}
uniformly for a.e. $r\in I$ and $\mu$ on bounded sets.

\indent Under the condition of $m\in M(I)$ and (\ref{hg1}), we shall show that
$\left(\mu_k^\nu(p),0\right)$ is a bifurcation point of (\ref{eg1}) and there are two distinct unbounded continua,
$\left(\mathcal{C}_{k}^\nu\right)^+$ and $\left(\mathcal{C}_{k}^\nu\right)^-$,
consisting of the bifurcation branch $\mathcal{C}_{k}^\nu$ from $\left(\mu_k^\nu(p), 0\right)$,
where $\mu_k^\nu(p)$ is the $k$-th positive or negative eigenvalue of
the linear problem corresponding to (\ref{eg1}), where $\nu\in\{+,-\}$.

Based on the unilateral global bifurcation result (see Theorem 3.2), we
investigate the existence of radial nodal solutions for the following $p$-Laplacian problem
\begin{equation}\label{eBf}
\left\{
\begin{array}{l}
-\text{div}\left(\varphi_p(\nabla u)\right)=\gamma m(x)f(u),\,\, \text{in\,}B,\\
u=0,\quad\quad\quad\quad\quad\quad\quad\quad\quad\quad\,\,\text{on}\, \partial B,
\end{array}
\right.
\end{equation}
where $f\in C(\mathbb{R})$, $\gamma$ is a parameter.

It is clear that the radial solutions of (\ref{eBf}) is
equivalent to the solutions of the following problem
\begin{equation}\label{erf1}
\left\{
\begin{array}{l}
\left(r^{N-1}\vert u'\vert^{p-2}u'\right)'+\gamma r^{N-1}m(r)f(u)=0,\,\, r\in I,\\
u'(0)=u(1)=0,
\end{array}
\right.
\end{equation}
where $r=\vert x\vert$ with $x\in B$.

It is well known that when $m(r)\equiv 1$ and $f(r, u) =\lambda\varphi_p(u)/\gamma$,
problem (\ref{erf1}) has a nontrivial
solution if and only if $\lambda$ is an eigenvalue of the following problem
\begin{equation}\label{eerp}
\left\{
\begin{array}{l}
\left(r^{N-1}\vert u'\vert^{p-2}u'\right)'+\lambda r^{N-1}\varphi_p(u)=0, \,\, r\in
I,\\
u'(0)=u(1)=0.
\end{array}
\right.
\end{equation}
In particular, when $\lambda=\lambda_k(p)$, there exist two solutions $u_k^+$ and
$u_k^-$, such that $u_k^+$
has exactly $k-1$ zeros in $I$ and is positive near 0, and $u_k^-$
has exactly $k-1$ zeros in $I$ and is negative near 0 (see [\ref{P}], Theorem 1.5.3).

When $p=2$, $N=1$ and $m(r)\geq 0$, Ma and Thompson [\ref{MT1}]
considered the interval of $\gamma$, in which there exist nodal solutions of (\ref{erf1})
under some suitable assumptions on $f$. The results of the above have been extended to the the case of weight
function changes its sign by Ma and Han [\ref{MH}]. The results they obtained
extended some well known theorems of the
existence of positive solutions for related problems [\ref{Er}, \ref{EW}, \ref{HW}] and sign-changing solutions [\ref{NT}].
For the case $p\not\equiv 2$ but $N=1$, $m(r)\geq 0$, Dai and Ma [\ref{DM}] proved the existence of nodal solutions for (\ref{erf1}).

However, few results on the existence of radial nodal solutions, even positive solutions,
have been established for $N$-dimensional $p$-Laplacian problem with sign-changing
weight $m(r)$ on the unit ball of $\mathbb{R}^N$. In this paper, we shall establish
a similar result to Ma and Thompson [\ref{MT1}]
for $N$-dimensional $p$-Laplacian problem with sign-changing weight. Problem with sign-changing
weight arises from the selection-migration model in population genetics. In
this model, $m(r)$ changes sign corresponding to the fact that an allele
$A_1$ holds an advantage over a rival allele $A_2$ at
same points and is at a disadvantage at others; the parameter $r$
corresponds to the reciprocal of diffusion, for detail, see [\ref{F}].
For the applications of nodal solutions, see Lazer and McKenna [\ref{LM}] and Kurth [\ref{K}].

In high dimensional general domain case, based on Dr\'{a}bek and Huang's results (note their results also valid for bounded smooth domain), we shall investigate the existence of one-sign solutions for the problem (\ref{eBf}) with
$1<p<N$ and the general smooth domain $\Omega\subset \mathbb{R}^N$ with $N\geq 2$, i.e.,
\begin{equation}\label{eBfg}
\left\{
\begin{array}{l}
-\text{div}\left(\varphi_p(\nabla u)\right)=\gamma m(x)f(u),\,\, \text{in\,}\Omega,\\
u=0,\quad\quad\quad\quad\quad\quad\quad\quad\quad\quad\,\,\text{on}\, \partial \Omega.
\end{array}
\right.
\end{equation}
By a solution of (\ref{eBfg}) we understand $u\in W_0^{1,p}(\Omega)$ satisfying (\ref{eBfg}) in the weak sense.

The rest of this paper is arranged as follows. In Section 2, we
establish the eigenvalue theory of second order $p$-Laplacian Dirichlet
boundary value problem in the radial case with sign-changing weight.
In Section 3, we establish the unilateral global bifurcation theory for (\ref{eg1}).
In Section 4, we prove the existence of nodal solutions for (\ref{erf1}).
In Section 5, we study the existence of one-sign solutions for (\ref{eBfg}).

\section{Some preliminaries}

\quad\, In [\ref{MYZ}], by Pr\"{u}fer transformation, Meng, Yan and Zhang established the
spectrum of one-dimensional $p$-Laplacian
with an indefinite integrable weight. When $m\equiv1$, using oscillation method,
Peral [\ref{P}] established the eigenvalue theory of $N$-dimensional
$p$-Laplacian on the unite ball. However, applying their methods to $N$-dimensional $p$-Laplacian
on the unit ball with indefinite weight is very difficult, even cann't be used. In [\ref{ACM}],
using variational method, Anane, Chakrone and Monssa established the spectrum of one-dimensional $p$-Laplacian
with an indefinite weight. While, we don't know whether or not the eigenvalue function
$\mu_k^\nu(p)$ is continuous with respect to $p$, which was obtained by Anane, Chakrone and Monssa.
In this Section, by similar method of Anane, Chakrone and Monssa's,
we can establish the eigenvalue theory of second order $p$-Laplacian Dirichlet
boundary value problem in the radial case with indefinite weight. It is well known the continuity of eigenvalues with respect to $p$ is very
important in the studying of the global bifurcation phenomena for
$p$-Laplacian problems, see [\ref{DPM}, \ref{DPEM}, \ref{LS}, \ref{P}]. In this Section,
we shall also show that $\mu_k^\nu(p)$ is continuous with respect to $p$. Moreover,
we also establish a key lemma which will be used in Section 4.
\\

\indent Applying the similar method to prove [\ref{ACM}, Theorem 1] with obvious changes, we can obtain the following:\\ \\
\textbf{Theorem 2.1.} \emph{Assume $m\in M(I)$. The eigenvalue problem
\begin{equation}\label{eer2}
\left\{
\begin{array}{l}
\left(r^{N-1}\vert u'\vert^{p-2}u'\right)'+\mu m(r)r^{N-1}\vert u\vert^{p-2}u=0,\,\, r\in I,\\
u'(0)=u(1)=0
\end{array}
\right.
\end{equation}
has two infinitely many simple real eigenvalues
\begin{equation}
0<\mu_{1}^+(p)<\mu_{2}^+(p)<\cdots<\mu_{k}^+(p)<
 \cdots,\lim_{k\rightarrow+\infty}\mu_{k}^+(p)=+\infty,\nonumber
\end{equation}
\begin{equation}
0>\mu_{1}^-(p)>\mu_{2}^-(p)>\cdots>\mu_{k}^-(p)>\cdots,
\lim_{k\rightarrow+\infty}\mu_{k}^-(p)=-\infty\nonumber
\end{equation}
and no other eigenvalues.
Moreover,}\\

\emph{1. Every eigenfunction corresponding to eigenvalue ${\mu}_k^\nu(p)$, has exactly $k-1$
zeros.}\\

\emph{2. For every $k$, ${\mu}_k^\nu(m)$ verifies the strict monotonicity property with
respect to the weight $m$.}
\\ \\
\textbf{Remark 2.1.} Using Gronwall inequality [\ref{E}], we can easily show that all zeros of eigenfunction corresponding to eigenvalue ${\mu}_k^\nu(p)$ is simple.
\\

\indent
We first show that the principle eigenvalue function $\mu_1^\nu:(1,+\infty)\rightarrow \mathbb{R}$ is continuous.
\\ \\
\textbf{Theorem 2.2.} \emph{The eigenvalue function $\mu_1^\nu:(1,+\infty)\rightarrow \mathbb{R}$ is continuous.}
\\ \\
\textbf{Proof.} The proof is similar to the proof of [\ref{DPM}], but we give a rough sketch of the proof for reader's
convenience. We only show that $\mu_1^+:(1,+\infty)\rightarrow \mathbb{R}$ is continuous since the case of $\mu_1^-$ is
similar. In the following proof, we shall shorten  $\mu_1^+$ to $\mu_1$.

From the variational characterization of $\mu_1(p)$ it follows that
\begin{equation}\label{cp1}
\mu_1(p)=\sup\left\{\mu>0\Big|\mu\int_Bm(x)\vert u\vert^p\,dx\leq\int_B\vert\nabla u\vert^p\,dx,
\,\,\text{for all\,\,}u\in C_{r,c}^\infty(B)\right\},
\end{equation}
where $C_{r,c}^\infty(B)=\left\{u\in C_{c}^\infty(B)\big| u\,\,\text{is radially symmetric}\right\}$.

Let $\{p_j\}_{j=1}^\infty$ be a sequence in $(1, +\infty)$ convergent to $p > 1$. We shall show
that
\begin{equation}\label{cp2}
\lim_{j\rightarrow+\infty}\mu_1(p_j)=\mu_1(p).
\end{equation}
To do this, let $u\in C_{r,c}^\infty(I)$. Then, from (\ref{cp1}),
\begin{equation}
\mu_1(p_j)\int_Bm(x)\vert u\vert^{p_j}\,dx\leq\int_B\vert\nabla u\vert^{p_j}\,dx.\nonumber
\end{equation}
On applying the Dominated Convergence Theorem we find
\begin{equation}\label{cp3}
\limsup_{j\rightarrow+\infty}\mu_1(p_j)\int_Bm(x)\vert u\vert^{p}\,dx\leq\int_B\vert\nabla u\vert^{p}\,dx.
\end{equation}
Relation (\ref{cp3}), the fact that $u$ is arbitrary and (\ref{cp1}) yield
\begin{equation}
\limsup_{j\rightarrow+\infty}\mu_1(p_j)\leq\mu_1(p).\nonumber
\end{equation}

Thus, to prove (\ref{cp2}) it suffices to show that
\begin{equation}\label{cp4}
\liminf_{j\rightarrow+\infty}\mu_1(p_j)\geq\mu_1(p).
\end{equation}
Let $\{p_k\}_{k=1}^\infty$ be a subsequence of $\{p_j\}_{j=1}^\infty$ such that
$\underset{k\rightarrow+\infty}\lim\mu_1(p_k)=\underset{j\rightarrow+\infty}\liminf\mu_1(p_j)$.

Let us fix $\varepsilon_0>0$ so that $p-\varepsilon_0>1$ and for each $0<\varepsilon<\varepsilon_0$,
$W_{r,0}^{1,p-\varepsilon}(B)$ is compactly embedded into $L_r^{p+\varepsilon}(B)$,
here $W_{r,0}^{1,p-\varepsilon}(B)=\left\{u\in W_{0}^{1,p-\varepsilon}(B)\big| u\,\,\text{is radially symmetric}\right\}$,
$L_r^{p+\varepsilon}(B)=\{u\in L_r^{p+\varepsilon}(B)\big| u\,\,\text{is radially symmetric}\}$.
For $k\in \mathbb{N}$, let us choose $u_k\in W_{r,0}^{1,p_k}(B)$ such that
\begin{equation}\label{cp5}
\int_B\vert\nabla u_k\vert^{p_k}\,dx=1
\end{equation}
and
\begin{equation}\label{cp6}
\int_B\vert\nabla u_k\vert^{p_k}\,dx=\mu_1(p_k)\int_Bm(x)\vert u_k\vert^{p_k}\,dx.
\end{equation}
For $0<\varepsilon<\varepsilon_0$, there exists $k_0\in \mathbb{N}$ such that
$p-\varepsilon<p_k<p+\varepsilon$ for any $k\geq k_0$. Thus, for $k\geq k_0$, (\ref{cp5}) and H\"{o}lder's inequality imply that
\begin{equation}\label{cp7}
\int_B\vert\nabla u_k\vert^{p-\varepsilon}\,dx\leq \vert B\vert^{\frac{p_k-p+\varepsilon}{p_k}},
\end{equation}
where $\vert B\vert$ denotes the measure of $B$.
This shows that $\{u_k\}_{k=k_0}^\infty$ is a bounded sequence in $W_{r,0}^{1,p-\varepsilon}(B)$. Passing to a
subsequence if necessary, we can assume that $u_k \rightharpoonup u$ in $W_{r,0}^{1,p-\varepsilon}(B)$ and hence
that $u_k \rightarrow u$ in $L_r^{p+\varepsilon}(B)$. Furthermore, $u\in L_r^{p}(B)$ and $u_k \rightarrow u$
in $L_r^{p_k}(B)$ for $k\geq k_0$.
It follows that
\begin{eqnarray}
\left\vert\int_B\vert u_k\vert^{p_k}\,dx-\int_B\vert u\vert^{p_k}\,dx\right\vert
&\leq&\int_B p_k\vert u+\theta u_k\vert^{p_k-1}\vert u_k-u\vert\,dx\nonumber\\
&\leq&(p+\varepsilon)\left(\int_B \vert u+\theta u_k\vert^{p_k}\,dx\right)^{\frac{p_k-1}{p_k}}\left(\int_B\vert u_k-u\vert^{p_k}\,dx\right)^{\frac{1}{p_k}}\nonumber\\
&\leq&(p+\varepsilon)\left(\Vert u\Vert_{p_k}+\Vert u_k\Vert_{p_k}\right)^{p_k-1}\left(\int_B\vert u_k-u\vert^{p_k}\,dx\right)^{\frac{1}{p_k}}\rightarrow0\nonumber
\end{eqnarray}
as $k\rightarrow +\infty$.
It is clear that
\begin{equation}
\int_B\vert u\vert^{p_k}\,dx-\int_B\vert u\vert^p\,dx\rightarrow 0\,\,\text{as}\,\,k\rightarrow +\infty.\nonumber
\end{equation}
Thus,
\begin{equation}
\int_B \vert u_k\vert^{p_k}\,dx\rightarrow\int_B\vert u\vert^{p}\,dx.\nonumber
\end{equation}
Similarly, we also can obtain that
\begin{equation}
\int_B m^+(x)\vert u_k\vert^{p_k}\,dx\rightarrow\int_Bm^+(x)\vert u\vert^{p}\,dx\nonumber
\end{equation}
and
\begin{equation}
\int_B m^-(x)\vert u_k\vert^{p_k}\,dx\rightarrow\int_Bm^-(x)\vert u\vert^{p}\,dx,\nonumber
\end{equation}
where $m^+(x)=\max\{m(x),0\}$, $m^-(x)=-\min\{m(x),0\}$.
Therefore,
\begin{eqnarray}\label{cp8}
\int_B m(x)\vert u_k\vert^{p_k}\,dx&=&\int_B m^+(x)\vert u_k\vert^{p_k}\,dx-\int_B m^-(x)\vert u_k\vert^{p_k}\,dx\nonumber\\
&\rightarrow&\int_Bm^+(x)\vert u\vert^{p}\,dx-\int_Bm^-(x)\vert u\vert^{p}\,dx\nonumber\\
&=&\int_Bm(x)\vert u\vert^{p}\,dx.
\end{eqnarray}
We note that (\ref{cp5}) and (\ref{cp6}) imply that
\begin{equation}\label{cp9}
\mu_1(p_k)\int_B m(x)\vert u_k\vert^{p_k}\,dx=1
\end{equation}
for all $k\in \mathbb{N}$. Thus letting $k$ go to $+\infty$ in (\ref{cp9}) and using (\ref{cp8}), we find
\begin{equation}\label{cp10}
\liminf_{j\rightarrow+\infty}\mu_1(p_k)\int_Bm(x)\vert u\vert^{p}\,dx=1.
\end{equation}
On the other hand, since $u_k\rightharpoonup u$ in $W_{r,0}^{1,p-\varepsilon}(B)$, from (\ref{cp7}) we obtain that
\begin{equation}
\Vert\nabla u\Vert_{p-\varepsilon}^{p-\varepsilon}\leq\liminf_{k\rightarrow+\infty}
\Vert\nabla u_k\Vert_{p-\varepsilon}^{p-\varepsilon}\leq \vert B\vert^{\frac{\varepsilon}{p}}.\nonumber
\end{equation}
Now, letting $\varepsilon\rightarrow 0^+$ and applying Fatou's Lemma we find
\begin{equation}\label{cp11}
\Vert\nabla u\Vert_p^p\leq 1.
\end{equation}
Hence $u\in W_r^{1,p}(B)$, here $W_r^{1,p}(B)$ denotes the radially symmetric subspace of $W^{1,p}(B)$.
We claim that actually $u\in W_{r,0}^{1,p}(B)$. Indeed, we know that
$u\in W_{r,0}^{1,p-\varepsilon}(B)$ for each $0<\varepsilon<\varepsilon_0$.
For $\phi\in C_{r,c}^\infty\left(\mathbb{R}^N\right)$ it is easy to see that
\begin{equation}
\left\vert \int_B u\frac{\partial\phi}{\partial x_i}\,dx\right
\vert\leq\Vert\nabla u\Vert_{p-\varepsilon}\Vert \phi\Vert_{(p-\varepsilon)'},\,\, i=1,\ldots,N.\nonumber
\end{equation}
Then, letting $\varepsilon\rightarrow 0^+$ we obtain that
\begin{equation}
\left\vert \int_B u\frac{\partial\phi}{\partial x_i}\,dx\right\vert\leq
\Vert\nabla u\Vert_{p}\Vert \phi\Vert_{(p)'},\,\, i=1,\ldots,N,\nonumber
\end{equation}
where $p'=p/(p-1)$.
Since $\phi$ is arbitrary, from Proposition IX-18 of [\ref{Bre1}] we find that $u\in W_{r,0}^{1,p}(B)$,
as desired.

Finally, combining (\ref{cp10}) and (\ref{cp11}) we obtain that
\begin{equation}
\liminf_{j\rightarrow+\infty}\mu_1(p_k)\int_Bm(x)
\vert u\vert^{p}\,dx\geq\int_B\vert\nabla u\vert^{p}\,dx.\nonumber
\end{equation}
This and the variational characterization of $\mu_1(p)$ imply (\ref{cp4}) and hence
(\ref{cp2}). This concludes the proof of the lemma.\qed\\

Using Remark 2.1, Theorem 2.1 and Theorem 2.2, we shall show that all eigenvalue
functions $\mu_k^\pm:(1,+\infty)\rightarrow \mathbb{R}$, $2\leq k\in \mathbb{N}$  are continuous.
\\ \\
\textbf{Theorem 2.3.} \emph{For each $2\leq k\in \mathbb{N}$, the eigenvalue
function $\mu_k^\nu:(1,+\infty)\rightarrow \mathbb{R}$ is continuous.}
\\ \\
\textbf{Proof.} Let $u_k^\nu$ be an eigenfunction corresponding to $\mu_k^\nu(p)$. By Theorem 2.1
and Remark 2.1, we know that $u$ has exactly $k-1$ simple zeros in $I$,
i.e., there exist $c_{k,1},\ldots,c_{k,k-1}\in I$ such that $u(c_{k,1})=\cdots=u(c_{k,k-1})=0$. For convenience,
we set $c_{k,0}=0$, $c_{k,k}=1$, $J_i=(c_{k,i-1},c_{k,i})$ and $B_i=\{x\in B\big|c_{k,i-1}<\vert x\vert<c_{k,i}\}$ for $i=1,\ldots, k$.
Let
$\mu_1^\nu\left(p,m/J_i,J_i\right)$ denote the first positive or negative eigenvalue
of the restriction of problem (\ref{eer2}) on $J_i$ for $i=1,\ldots, k$.

We note that Lemma 3 of [\ref{ACM}] also holds for (\ref{eer2}). It follows that $\mu_k^\nu(p)=\mu_1^\nu\left(p,m/J_i,J_i\right)$
for $i=1,\ldots,k$. Using similar proof as Theorem 2.2, we can show
that $\mu_1^\nu\left(p,m/J_i,J_i\right)$ is continuous with respect to $p$ for
$i=1,\ldots, k$. Therefore, $\mu_k^\nu(p)$ is also continuous with respect to $p$.  \qed
\\

\indent Finally, we give a key lemma that will be used in Section 4. Firstly, as an immediate consequence of Lemma 4.1 of [\ref{DPM}], we obtain the following Sturm type comparison theorem.\\ \\
\textbf{Lemma 2.1.} \emph{Let $b_2(r)>b_1(r)>0$ for $r\in(0,1)$ and $b_i(r)\in L^\infty(0,1)$, $i=1,2$. Also let $u_1$, $u_2$
are solutions of
\begin{equation}
\left(r^{N-1}\varphi_p(u')\right)'+b_i(r)r^{N-1}\varphi_p(u)=0,\,\, i=1,2,\nonumber
\end{equation}
respectively. If $u_1$ has $k$ zeros in $(0,1)$, then $u_2$ has at least $k+1$ zeros in $(0,1)$.}
\\ \\
\indent Let
\begin{equation}
I^+:=\left\{r\in \overline{I}\,|\, m(r)>0\right\}, \ \ \ I^-:=\left\{r\in \overline{I}\,|\, m(r)<0\right\}.\nonumber
\end{equation}
\noindent\textbf{Lemma 2.2.} \emph{Assume $m\in M(I)$. Let $\widehat{I}=[a,b]$ be such
that $\widehat{I}\subset I^+$ and
\begin{equation}
\text{meas}\,\widehat{I}>0.
\nonumber
\end{equation}
Let $g_n:\overline{I}\to (0, +\infty)$ be continuous function and such that
\begin{equation}
\lim_{n\to +\infty} g_n(r)=+\infty\ \ \text{uniformly on}\
\widehat{I}.\nonumber
\end{equation}
Let $y_n\in E$ be a solution of the equation
\begin{equation}
\left(r^{N-1}\varphi_p(y_n')\right)'+ r^{N-1}m(r)g_n(r)\varphi_p(y_n)=0, \,\, r\in (0, 1).\nonumber
\end{equation}
Then the number of zeros of $y_n|_{\widehat{I}}$ goes to infinity as $n\to +\infty$.}
\\ \\
\noindent\textbf{Proof.} After taking a subsequence if
necessary, we may assume that
\begin{equation}
m(r)g_{n_j}(r)\geq \lambda_j, \ \  r\in \widehat{I},\nonumber
\end{equation}
as $j\to +\infty$, where $\lambda_j$ is the $j$-th eigenvalue of the following problem
\begin{equation}
\left\{
\begin{array}{l}
\left(r^{N-1}\varphi_p(u'(r))\right)'+\lambda r^{N-1}\varphi_p(u(r))=0,\,\,r\in I,\\
u'(0)=u(1)=0.
\end{array}
\right.\nonumber
\end{equation}
Let $\varphi_j$ be the corresponding eigenvalue of $\lambda_j$.
It is easy to check that the number of zeros of $\varphi_j\big|_{\widehat{I}}$ goes to infinity as $j\to +\infty$.
By Lemma 2.1, one has that the number of zeros of $y_n|_{\widehat{I}}$ goes to infinity as $n\to +\infty$.
It follows the desired results.\qed
\section{Unilateral global bifurcation phenomena for (\ref{eg1})}
\quad\, If $m(r)\equiv1$, Del Pino and Elgueta [\ref{DPEM}] established
the global bifurcation theory for one dimensional $p$-Laplacian
eigenvalue problem. Peral [\ref{P}] got the global bifurcation
theory for $p$-Laplacian eigenvalue problem on the unite ball. In
[\ref{DPM}], Del Pino and Man\'{a}sevich obtained the global
bifurcation from the principle eigenvalue for $p$-Laplacian
eigenvalue problem on the general domain. If $m(r)\geq 0$ and is
singular at $r=0$ or $r=1$, Lee and Sim [\ref{LS}] also
established the bifurcation theory for one dimensional $p$-Laplacian
eigenvalue problem. However, if $m(r)$ changes sign, there are a few
paper involving in the bifurcation theory for $p$-Laplacian
eigenvalue problem. In this Section, we shall study the unilateral global bifurcation phenomena
for $N$-dimensional $p$-Laplacian eigenvalue problem with sign-changing weight in the radial case.

Let $Y=L^1(0,1)$ with its usual normal $||\cdot||_{L^1}$
and $E = \left\{u \in C^1(\overline{I})\big| u'(0) = u(1) = 0\right\}$ with the norm
\begin{equation}
\Vert u\Vert=\max_{r\in\overline{I}}\vert u(r)\vert+\max_{r\in\overline{I}}\vert u'(r)\vert.\nonumber
\end{equation}

Considering the following auxiliary problem
\begin{equation}\label{erh}
\left\{
\begin{array}{l}
-\left(r^{N-1}\vert u'\vert^{p-2}u'\right)'=r^{N-1}h(r),\,\, \text{ a.e.}\,\, r\in I,\\
u'(0)=u(1)=0
\end{array}
\right.
\end{equation}
for a given $h\in Y$.
By a solution of problem (\ref{erh}), we understand a function $u\in E$ with $r^{N-1}\varphi_p(u')$ absolutely
continuous which satisfies (\ref{erh}).

We have known that for every given $h\in Y$ there is a unique solution
$u$ to the problem (\ref{erh}) (see [\ref{DPM}]). Let $G_p(h)$ denote the unique solution to (\ref{erh})
for a given $h\in Y$.
It is well known that $G_p:Y\rightarrow E$ is continuous and compact (see [\ref{DPM}, \ref{P}]).
\\

\indent Define $T_\mu^p(u)=G_p\left(\mu m(r)\varphi_p(u(r))\right)$. Let $\Psi_{p,\mu}$ be defined in $E$ by
\begin{equation}
\Psi_{p,\mu}(u)=u-T_\mu^p(u),\nonumber
\end{equation}
where $\mu$ is a positive parameter.
It is no difficult to show that $\Psi_{p,\mu}$ is a nonlinear
compact perturbation of the identity. Thus the Leray-Schauder degree
$\deg\left(\Psi_{p,\mu}, B_r(0),0\right)$ is well-defined for
arbitrary $r$-ball $B_r(0)$ and $\mu\neq \mu_k^\nu$. \\

\indent Firstly, we
can compute $\deg\left(\Psi_{2,\mu}, B_r(0),0\right)$ for any $r>0$ as follows.
\\ \\
\noindent\textbf{Lemma 3.1.} \emph{For $r>0$, we have
$$\deg\left(\Psi_{2,\mu}, B_r(0),0\right)
=\left\{\aligned 1, \ \ \ \qquad \ \ & \text{if}\ \mu\in
\left(\mu_1^-(2),
\mu_1^+(2)\right),\\
(-1)^k, \ \ \ \ \ \ & \text{if}\ \mu\in \left(\mu_k^+(2), \mu_{k+1}^+(2)\right),\ k\in \mathbb{N},\\
(-1)^k, \ \ \ \ \ \ & \text{if}\ \mu\in \left(\mu_{k+1}^-(2), \mu_{k}^-(2)\right),\ k\in \mathbb{N}.\\
\endaligned
\right.
$$
}
\\ \\
\noindent\textbf{Proof.}~ We divide the proof into two cases.

{\it Case 1.} $\mu\geq 0$.

Since $G_2$ is compact and linear, by [\ref{De}, Theorem 8.10] and
Theorem 2.1 with $p=2$,
\begin{equation}
\deg\left(\Psi_{2,\mu}, B_r(0),0\right)=(-1)^{m(\mu)},\nonumber
\end{equation}
where $m(\mu)$ is the sum of algebraic multiplicity of the
eigenvalues $\mu$ of (\ref{eer2}) satisfying $\mu^{-1}\mu_k^+<1$.
If $\mu\in [0, \mu^+_1(2))$, then there are no such $\mu$ at
all, then
\begin{equation}
\deg\left(\Psi_{2,\mu}, B_r(0),0\right)=(-1)^{m(\mu)}=(-1)^0=1.\nonumber
\end{equation}
If $\mu\in \left(\mu_k^+(2), \mu_{k+1}^+(2)\right)$ for some $k\in
\mathbb{N}$, then
\begin{equation}
\left(\mu^+_j(2)\right)^{-1}\mu>1, \ \  j\in \{1, \cdots, k\}.\nonumber
\end{equation}
This together with Theorem 2.1 implies
\begin{equation}
\deg\left(\Psi_{2,\mu}, B_r(0),0\right)=(-1)^k.\nonumber
\end{equation}

{\it Case 2.} $\mu<0$.

In this case, we consider a new sign-changing eigenvalue problem
$$
\left\{
\begin{array}{l}
\left(r^{N-1}u'\right)'+\hat\mu \hat m(r)r^{N-1}u=0,\,\,r\in I,\\
u'(0)=u(1)=0,
\end{array}
\right.
$$
where $\hat\mu=-\mu$, $\hat m(r)=-m(r)$. It is easy to check
that
\begin{equation}
\hat\mu_k^+(2)=-\mu_k^-(2), \ \  k\in
\mathbb{N}.\nonumber
\end{equation}
Thus, we may use the result obtained in \emph{Case 1} to deduce the desired
result. \qed\\

\indent As far as the general $p$ is concerned, we can compute the extension of the Leray-Schauder degree defined in [\ref{BP}]
by the deformation along $p$.\\ \\
\noindent\textbf{Lemma 3.2.} (i) \emph{Let
$\left\{\mu_{k}^{+}(p)\right\}_{k\in \mathbb{N}}$ be the sequence of
positive eigenvalues of (\ref{eer2}). Let $\mu$ be a constant with
$\mu\neq\mu_{{k}}^+(p)$ for all $k\in \mathbb{N}$. Then for
arbitrary $r>0$,
\begin{equation}
\deg \left(\Psi_{p,\mu}, B_r(0),0 \right)=(-1)^{\beta},\nonumber
\end{equation}
where $\beta$ is the number of eigenvalues $\mu_k^+(p)$ of problem
(\ref{eer2}) less than $\mu$.}

(ii) \emph{Let $\{\mu_{k}^{-}(p)\}_{k\in \mathbb{N}}$ be the
sequence of negative eigenvalues of (\ref{eer2}). Consider
$\mu\neq\mu_{{k}}^-(p)$, $k\in \mathbb{N}$, then
\begin{equation}
\deg\left(\Psi_{p,\mu}, B_r(0),0\right)=(-1)^{\beta},\,\, \forall r>0,\nonumber
\end{equation}
where $\beta$ is the number of eigenvalues $\mu_k^-(p)$ of problem
(\ref{eer2}) larger than $\mu$.}
\\ \\
\noindent\textbf{Proof.} We shall only prove the case $\mu>\mu_1^+(p)$
since the proof for the other cases are similar. We
also only give the proof for the case $p>2$. Proof for the case
$1<p<2$ is similar. Assume that $\mu_{k}^+(p) < \mu<
\mu_{{k+1}}^+(p)$ for some $k\in\mathbb{N}$. Since the eigenvalues
depend continuously on $p$, there exists a continuous function
$\chi:[2,p]\rightarrow\mathbb{R}$ and $q\in [2,p]$ such that
$\mu_{k}^+(q) < \chi(q) < \mu_{{k+1}}^+(q)$ and $\lambda=\chi(p)$.
Define
\begin{equation}
\Upsilon(q,u)=u-G_q\left(\chi(q)m(r)\varphi_q(u)\right).\nonumber
\end{equation}
It is easy to show that $\Upsilon(q,u)$ is a compact perturbation of the identity such that for all $u\neq 0$,
by definition of $\chi(q)$,
$\Upsilon(q,u)\neq0$, for all $q\in [2,p]$. Hence the invariance of the degree under homo-topology and Lemma 3.1 imply
\begin{equation}
\deg\left(\Psi_{p,\mu}, B_r(0),0\right)=\deg\left(\Psi_{2,\mu},
B_r(0),0\right)=(-1)^{k}.\nonumber
\end{equation} \qed

\indent Define the Nemitskii operator $H:\mathbb{R}\times E\rightarrow Y$
by
$$
H(\mu,u)(r):=\mu m(r)\varphi_p(u(r))+g(r,u(r);\mu).\nonumber
$$
Then it is clear that $H$ is continuous (compact) operator and problem (\ref{eg1}) can be equivalently written as
$$
u=G_p\circ H(\mu,u):=F(\mu,u).
$$
$F$ is completely continuous in $\mathbb{R}\times E\rightarrow E$ and
$F(\mu,0)=0$, $\forall \mu\in \mathbb{R}$.\\

Using the similar method to prove [\ref{DM}, Theorem 2.1] with obvious changes, we may obtain the
following result.
\\ \\
\textbf{Theorem 3.1.} \emph{Assume (\ref{hg1}) holds and $m\in M(I)$, then from each $\left(\mu_k^\nu,0\right)$ it
bifurcates an unbounded continuum $\mathcal{C}_k^\nu$ of solutions to problem (\ref{eg1}),
with exactly $k-1$ simple zeros, where $\mu_k^\nu$ is the eigenvalue of problem (\ref{eer2}).}\\

In what follows, we use the terminology of Rabinowitz [\ref{R3}]. Let $S_k^+$ denote the set of
functions in $E$ which have exactly $k-1$ interior
nodal (i.e. non-degenerate zeros) in $I$ and are positive near $t=0$, and set $S_k^-=-S_k^+$, and
$S_k =S_k^+\cup S_k^-$. They are disjoint and open in $E$. The following global bifurcation result is a generalization of Theorem 3.2 of [\ref{DM}]. The essential idea is similar to the proof of
Theorem 3.2 of [\ref{DM}]\\ \\
\textbf{Theorem 3.2.} \emph{Assume (\ref{hg1}) holds and $m\in M(I)$, then there are two distinct unbounded continua,
$\left(\mathcal{C}_{k}^\nu\right)^+$ and $\left(\mathcal{C}_{k}^\nu\right)^-$,
consisting of the bifurcation branch $\mathcal{C}_{k}^\nu$. Moreover, for $\sigma\in\{+,-\}$, we have
\begin{equation}
\left(\mathcal{C}_{k}^\nu\right)^\sigma\subset \left(\{\left(\mu_k^\nu,0\right)\}\cup\left(\mathbb{R}\times S_k^\sigma\right)\right).\nonumber
\end{equation}}

\section{Existence of nodal solutions of (\ref{erf1})}

\quad\, In this Section, we shall investigate the existence and multiplicity of nodal
solutions to the problem (\ref{erf1}) under the linear growth condition on $f$.\\

\indent Firstly, we suppose that\\

\emph{(${H}_1$) $f\in C(\mathbb{R},\mathbb{R})$ with
$f(s)s>0$ for $s\neq0$;}

\emph{($H_2$) there exists $f_0\in (0, +\infty)$ such that}
\begin{equation}
f_0=\lim_{\vert s\vert\rightarrow0}\frac{f(s)}{\varphi_p(s)};\nonumber
\end{equation}

\emph{($H_3$) there exists $f_\infty \in (0,+\infty)$ such
that}
\begin{equation}
f_\infty=\lim_{\vert
s\vert\rightarrow+\infty}\frac{f(s)}{\varphi_p(s)}.\nonumber
\end{equation}

\indent Let $\mu_k^\pm$ be the $k$-th positive or negative eigenvalue of (\ref{eer2}).
Applying Theorem 3.2, we shall establish the existence of nodal solutions of (\ref{erf1}) follows.\\ \\
\textbf{Theorem 4.1.} \emph{Assume ($H_1$), ($H_2$) and ($H_3$) hold and $m\in M(I)$.
Assume that for some $k \in \mathbb{N}$, either
\begin{equation}
\gamma\in\left(\frac{\mu_k^+(p)}{f_\infty},\frac{\mu_k^+(p)}{f_0}\right)
\cup\left(\frac{\mu_k^-(p)}{f_0},\frac{\mu_k^-(p)}{f_\infty}\right)\nonumber
\end{equation}
or
\begin{equation}
\gamma\in\left(\frac{\mu_k^+(p)}{f_0},\frac{\mu_k^+(p)}{f_\infty}\right)
\cup\left(\frac{\mu_k^-(p)}{f_\infty},\frac{\mu_k^-(p)}{f_0}\right).\nonumber
\end{equation}
Then (\ref{erf1}) has two solutions $u_k^+$ and $u_k^-$ such
that $u_k^+$ has exactly $k-1$ zeros in $I$ and is positive near 0,
and $u_k^-$ has exactly $k-1$ zeros in $I$ and is negative near 0.}
\\ \\
\textbf{Proof.} We only prove the case of $\gamma>0$. The case of
$\gamma<0$ is similar. Consider the problem
\begin{equation}\label{ff}
\left\{
\begin{array}{l}
\left(r^{N-1}\varphi_p(u')\right)'+\mu\gamma r^{N-1}{m}(r)f(u)=0,\,\, r\in I,\\
u'(0)=u(1)=0.
\end{array}
\right.
\end{equation}
Let $\zeta\in C(\mathbb{R})$ be such that $f(u)=f_0\varphi_p(u)+\zeta(u)$
with $\lim_{\vert u\vert\rightarrow0}\zeta(u)/\varphi_p(u)=0.$
Hence, the condition
(\ref{hg1}) holds. Using Theorem 3.2, we have that
there are two distinct unbounded continua, $\left(\mathcal{C}_k^\nu\right)^+$ and $\left(\mathcal{C}_k^\nu\right)^-$,
consisting of the bifurcation branch $\mathcal{C}_k^\nu$ from $\left(\mu_k^\nu(p)/\gamma, 0\right)$
, such that
\begin{equation}
\left(\mathcal{C}_{k}^\nu\right)^\sigma\subset \left(\left\{\left(\mu_k^\nu,0\right)\right\}\cup\left(\mathbb{R}\times S_k^\sigma\right)\right).\nonumber
\end{equation}

It is clear that any solution of (\ref{ff}) of the form $(1, u)$ yields a
solutions $u$ of (\ref{erf1}). We shall show that $\left(C_k^+\right)^\sigma$ crosses
the hyperplane $\{1\}\times E$ in $\mathbb{R}\times E$. To  this
end, it will be enough to show that $\left(C_k^+\right)^\sigma$ joins
$\left(\mu_k^+(p)/\gamma f_0, 0\right)$ to
$\left(\mu_k^+(p)/\gamma f_\infty, +\infty\right)$. Let
$\left(\eta_n, y_n\right) \in \left(C_k^+\right)^\sigma$ satisfy $\eta_n+\Vert y_n\Vert\rightarrow+\infty.$
We note that $\eta_n >0$ for all $n \in \mathbb{N}$ since (0,0) is
the only solution of (\ref{ff}) for $\mu = 0$ and
$\left(C_k^+\right)^\sigma\cap\left(\{0\}\times E\right)=\emptyset$.

\emph{Case 1}:
$\mu_k^+(p)/f_\infty<\gamma<\mu_k^+(p)/f_0$.
 In this case, we only need to show that
\begin{equation}
\left(\frac{\mu_k^+(p)}{\gamma f_\infty},\frac{\mu_k^+(p)}{\gamma
f_0}\right) \subseteq\left\{\mu\in\mathbb{R}\big|(\mu,u)\in \left(C_k^+\right)^\sigma\right\}.\nonumber
\end{equation}
We divide the proof into two steps.

\emph{Step 1}: We show that if there exists a constant $M>0$ such that $\eta_n\subset(0,M]$ for $n\in \mathbb{N}$ large enough, then $(C_k^+)^\sigma$ joins
$\left(\mu_k^+(p)/\gamma f_0, 0\right)$ to
$\left(\mu_k^+(p)/\gamma f_\infty, +\infty\right)$.

In this case it follows that $\Vert y_n\Vert\rightarrow+\infty.$
Let $\xi\in C(\mathbb{R})$ be such that $f(u)=f_\infty \varphi_p(u)+\xi(u).$
Then $\lim_{\vert u\vert\rightarrow+\infty}\frac{\xi(u)}{\varphi_p(u)}=0.$
Let $\widetilde{\xi}(u)=\max_{0\leq \vert s\vert\leq u}\vert \xi(s)\vert.$
Then $\widetilde{\xi}$ is nondecreasing and
\begin{equation}\label{eu0}
\lim_{u\rightarrow +\infty}\frac{\widetilde{\xi}(u)}{\vert
u\vert^{p-1}}=0.
\end{equation}

We divide the equation
\begin{equation}
\left(r^{N-1}\varphi_p(y_n')\right)'-\mu_n\gamma r^{N-1}m(r)f_\infty\varphi_p(y_n)=\mu_n\gamma r^{N-1}m(r)\xi(y_n)\nonumber
\end{equation}
by $\Vert y_n\Vert$ and set $\overline{y}_n = y_n/\Vert y_n\Vert$. Since $\overline{y}_n$ is bounded in $E$,
after taking a subsequence if
necessary, we have that $\overline{y}_n \rightharpoonup \overline{y}$ for some $\overline{y} \in E$. Moreover, from
(\ref{eu0}) and the fact that $\widetilde{\xi}$ is nondecreasing, we have that
\begin{equation}
\lim_{n\rightarrow+\infty}\frac{ \xi(y_n(r))}{\Vert y_n\Vert^{p-1}}=0\nonumber
\end{equation}
since
\begin{equation}
\frac{ \xi(y_n(r))}{\Vert y_n\Vert^{p-1}}\leq\frac{ \widetilde{\xi}(\vert y_n(r)\vert)}{\Vert y_n\Vert^{p-1}}
\leq\frac{ \widetilde{\xi}(\Vert y_n(r)\Vert_\infty)}{\Vert y_n\Vert^{p-1}}
\leq\frac{ \widetilde{\xi}(\Vert y_n(r)\Vert)}{\Vert y_n\Vert^{p-1}}.\nonumber
\end{equation}
By the continuity and compactness of $G_p$, it follows that
\begin{equation}
-\left(r^{N-1}\varphi_p(\overline{y}')\right)'=\overline{\mu}\gamma r^{N-1}m(r)f_\infty\varphi_p(\overline{y}),\nonumber
\end{equation}
where
$\overline{\mu}=\underset{n\rightarrow+\infty}\lim\mu_n$, again
choosing a subsequence and relabeling if necessary.

We claim that $\overline{y}\in \left(C_k^+\right)^\sigma.$

It is clear that $\Vert \overline{y}\Vert=1$ and $\overline{y}\in \overline{\left(\mathcal{C}_{k}^+\right)^\sigma}\subseteq
\left(\mathcal{C}_{k}^+\right)^\sigma$ since $\left(\mathcal{C}_{k}^+\right)^\sigma$ is closed in $\mathbb{R}\times E$.
Therefore, by Theorem 2.1, $\overline{\mu}\gamma f_\infty=\mu_k^+(p)$, so
that $\overline{\mu}=\mu_k/\gamma f_\infty.$
Therefore $(C_k^+)^\sigma$ joins $\left(\mu_k^+(p)/\gamma
f_0, 0\right)$ to $\left(\mu_k^+(p)/\gamma f_\infty,
+\infty\right)$.

\emph{Step 2}: We show that there exists a constant $M$ such that $\mu_n
\in(0,M]$ for $n\in \mathbb{N}$ large enough.

On the contrary, we suppose that $\lim_{n\rightarrow +\infty}\mu_n=+\infty.$
Since $\left(\eta_n, y_n\right) \in (C_k^+)^\sigma$, it follows that
\begin{equation}
\left(r^{N-1}\varphi(y_n')\right)'+\gamma\eta_n
r^{N-1}m(r)f(y_n)=0.\nonumber
\end{equation}
Let
\begin{equation}
0<\tau(1,n)<\cdots<\tau(k,n)=1\nonumber
\end{equation}
be the zeros of $y_n$ in $\overline{I}$. Then, after taking a subsequence
if necessary,
\begin{equation}
\lim_{n\rightarrow +\infty}\tau(l,n):=\tau(l,\infty), \qquad l\in
\{1, \cdots, k-1 \}.\nonumber
\end{equation}
It follows that either there exists at least one $l_0\in\{1, \cdots, k-1 \}$ such that
\begin{equation}
\tau(l_0,\infty)<\tau(l_0+1,\infty)\,\,\text{or}\,\,\tau(1,\infty)=1.\nonumber
\end{equation}
Notice that Lemma 2.2 and the fact $y_n$ has exactly $k-1$ simple
zeros in $\overline{I}$ yield
\begin{equation}
\left\{\left[\cup^{k-1}_{l=1}\big(\tau(l,\infty),\tau(l+1,\infty)\big)\right]\cup(0,\tau(1,\infty))\right\}\cap I^+=\emptyset,\nonumber
\end{equation}
which implies that
\begin{equation}
\left\{\left[\cup^{k-1}_{l=1}\big(\tau(l,\infty),\tau(l+1,\infty)\big)\right]\cup(0,\tau(1,\infty))\right\}\subseteq (I\setminus I^+).\nonumber
\end{equation}
Therefore,
\begin{equation}
\text{meas}\,(I\setminus I^+)\geq\text{meas}\,\left\{\left[\cup^{k-1}_{l=1}\big(\tau(l,\infty),\tau(l+1,\infty)\big)\right]\cup(0,\tau(1,\infty))\right\}=1.\nonumber
\end{equation}
However, this contradicts ($H_2$): $0<\text{meas}\,(I\setminus I^+)<1$.

\emph{Case 2}:
$\mu_k^+(p)/f_0<\gamma<\mu_k^+(p)/f_\infty$. In this
case, we have that
\begin{equation}
\frac{\mu_k^+(p)}{\gamma f_0}<1<\frac{\mu_k^+(p)}{\gamma
f_\infty}.\nonumber
\end{equation}
Assume that $\left(\eta_n, y_n\right) \in \left(C_k^+\right)^\sigma$ is such that $\lim_{n\rightarrow+\infty}\left(\eta_n+\Vert y_n\Vert\right)=+\infty.$
In view of \emph{Step 2} of \emph{Case 1}, we have known that there exists $M>0$, such that
for $n \in \mathbb{N}$ sufficiently large, $\eta_n\in (0,M].$
Applying the same method used in \emph{Step 1} of \emph{Case 1}, after
taking a subsequence and relabeling if necessary, it follows that
\begin{equation}
(\eta_n,y_n)\rightarrow\left(\frac{\mu_k^+(p)}{\gamma
f_\infty},+\infty\right) \,\, \text{as}\,\, n\rightarrow+\infty.\nonumber
\end{equation}
Thus, $(C_k^+)^\sigma$ joins $\left(\mu_k^+(p)/\gamma
f_0,0\right)$ to $\left(\mu_k^+(p)/\gamma
f_\infty,+\infty\right)$.\qed
\\

\indent Using the similar proof with the proof Theorem 4.1, we can obtain the more general results as follows.
\\ \\
\textbf{Theorem 4.2.} \emph{Assume ($H_1$), ($H_2$) and ($H_3$) hold and $m\in M(I)$.
Assume that for some $k, n \in \mathbb{N}$ with $k\leq n$, either
\begin{equation}
\gamma\in\left(\frac{\mu_n^+(p)}{f_\infty},\frac{\mu_k^+(p)}{f_0}\right)
\cup\left(\frac{\mu_k^-(p)}{f_0},\frac{\mu_n^-(p)}{f_\infty}\right)\nonumber
\end{equation}
or
\begin{equation}
\gamma\in\left(\frac{\mu_n^+(p)}{f_0},\frac{\mu_k^+(p)}{f_\infty}\right)
\cup\left(\frac{\mu_k^-(p)}{f_\infty},\frac{\mu_n^-(p)}{f_0}\right).\nonumber
\end{equation}
Then (\ref{erf1}) has $n-k+1$ pairs solutions $u_j^+$ and $u_j^-$ for $j\in\{k,\cdots,n\}$ such
that $u_j^+$ has exactly $j-1$ zero in $I$ and is positive near 0,
and $u_j^-$ has exactly $j-1$ zero in $I$ and is negative near 0.}
\\ \\
\textbf{Remark 4.1.} We would like to point out that Theorem 1.1 of [\ref{MT1}] is the corollary of Theorem 5.1
even in the case of $p=2$ and $N=1$.
\\ \\
\textbf{Remark 4.2.} We also note that Theorem 4.1 and Theorem 4.2 is valid for the problems on annular domain because
it can be convert the equivalent one-dimensional problems.\\ \\
\textbf{An open problem.} When $m\geq 0$, using Lemma 2.1, we can easily get that
(\ref{erf1}) has no nontrivial solution if $\gamma mf(u)/u$ not cross any eigenvalue of (\ref{eerp}). Therefore,
we conjecture that (\ref{erf1}) has no nontrivial solution if
\begin{equation}\label{ns}
\mu_k^+(p)<\frac{f(s)}{\varphi_p(s)}<\mu_{k+1}^+(p)\,\,\text{or\,\,} \mu_k^-(p)>-\frac{f(s)}{\varphi_p(s)}>\mu_{k+1}^-(p)\,\,\text{for}\,\, s\neq 0.\nonumber
\end{equation}

\section{One-sign solutions for (\ref{eBfg})}
\quad\, In this Section, based on the bifurcation result of Dr\'{a}bek and Huang [\ref{DH}], we
shall study the existence of one-sign solutions for
problem (\ref{eBfg}). From now on, for simplicity, we write $X:=W_0^{1,p}(\Omega)$.\\

The main results of this section are the following:\\ \\
\textbf{Theorem 5.1.} \emph{Let ($H_1$), ($H_2$) and ($H_3$) hold, and $m\in M(\Omega)$.
Assume that either
\begin{equation}
\gamma\in\left(\frac{\mu_1^+(p)}{f_\infty},\frac{\mu_1^+(p)}{f_0}\right)
\cup\left(\frac{\mu_1^-(p)}{f_0},\frac{\mu_1^-(p)}{f_\infty}\right)\nonumber
\end{equation}
or
\begin{equation}
\gamma\in\left(\frac{\mu_1^+(p)}{f_0},\frac{\mu_1^+(p)}{f_\infty}\right)
\cup\left(\frac{\mu_1^-(p)}{f_\infty},\frac{\mu_1^-(p)}{f_0}\right).\nonumber
\end{equation}
then problem (\ref{eBfg}) possesses at least a positive and a negative solution.}\\
\\
\textbf{Remark 5.1.} By the
$C^{1,\alpha}$ ($0<\alpha<1$) regularity results for quasilinear elliptic equations with $p$-growth condition [\ref{Li}], $u\in C^{1,\alpha}(\overline{\Omega})$ for any solution $u$ of (\ref{eBfg}) since $f$ is continuous and subcritical.
 \\

In order to prove Theorem 5.1, we consider the following eigenvalue problem
\begin{equation}\label{hdg1}
\left\{
\begin{array}{l}
-\text{div}\left(\varphi_p(\nabla u)\right)=\mu\gamma m(x)f(u),\,\, \text{in\,}\Omega,\\
u=0,\quad\quad\quad\quad\quad\quad\quad\quad\quad\quad\,\,\,\,\,\,\text{on}\, \partial \Omega,
\end{array}
\right.
\end{equation}
where $\mu$ is a parameter. Let $\zeta\in C(\mathbb{R})$ be such that
\begin{equation}
f(u)=f_0\varphi_p(u)+\zeta(u)\nonumber
\end{equation}
with $\lim_{\vert u\vert\rightarrow0}\zeta(u)/\varphi_p(u)=0$. Let us consider
\begin{equation}\label{hdg2}
\left\{
\begin{array}{l}
-\text{div}\left(\varphi_p(\nabla u)\right)=\mu\gamma m(x)f_0\varphi_p(u)+\mu\gamma m(x)\zeta(u),\,\, \text{in\,}\Omega,\\
u=0,\quad\quad\quad\quad\quad\quad\quad\quad\quad\quad\quad\quad\quad\quad\quad\quad\quad\quad\,\,\,\,\,\text{on}\, \partial \Omega,
\end{array}
\right.
\end{equation}
as a bifurcation problem from the trivial solution $u\equiv 0$.

Let
\begin{equation}
\mathbb{S}^+=\left\{u\in C^{1,\alpha}(\overline{\Omega})\big|u(x)> 0, \text{\,for all\,\,}x\in\Omega\right\}\,\,\text{and}\,\,\mathbb{S}^-=\left\{u\in C^{1,\alpha}(\overline{\Omega})\big|u(x)< 0, \text{\,for all\,\,}x\in\Omega\right\}.\nonumber
\end{equation}
Applying Theorem 4.4 and 4.5 of [\ref{DH}] to (\ref{hdg2}), we can obtain the following unilateral global bifurcation result, which plays a fundamental role in our study.\\ \\
\textbf{Lemma 5.1.} \emph{Let $\nu\in\{+,-\}$. There are two distinct unbounded continua, $\mathcal{C}_+^\nu$ and $\mathcal{C}_-^\nu$,
consisting of the bifurcation branch $\mathcal{C}^\nu$ from ($\mu_k^\nu(p), 0$). Moreover, for $\sigma\in\{+,-\}$, we have
\begin{equation}
\mathcal{C}_\sigma^\nu\subset \left(\{(\mu_1(p),0)\}\cup(\mathbb{R}\times \mathbb{S}^\sigma)\right).\nonumber
\end{equation}}
\indent We use Lemma 6.1 to prove the main results of this section.
\\ \\
\textbf{Proof of Theorem 5.1.}
Since the proof is similar to that of Theorem 4.1, we only give a rough sketch of the proof. We only prove the case of $\gamma>0$. The case of
$\gamma<0$ is similar. It is clear that any solution of (\ref{hdg1}) of the form $(1, u)$ yields a
solution $u$ of (\ref{eBfg}). We shall show $\mathcal{C}_\sigma^+$ crosses the hyperplane $\{1\}\times X$ in $\mathbb{R}\times X$.
To  this end, it will be enough to show that $\mathcal{C}_\sigma^+$ joins
$\left(\mu_1^+(p)/\gamma f_0, 0\right)$ to
$\left(\mu_1^+(p)/\gamma f_\infty, +\infty\right)$.

Let $(\mu_n, y_n) \in \mathcal{C}_\sigma^+$ where $y_n\not\equiv 0$ satisfies $\mu_n+\Vert y_n\Vert_X\rightarrow+\infty.$
We note that $\mu_n >0$ for all $n \in \mathbb{N}$ since (0,0) is the only solution of (\ref{hdg1}) for $\mu = 0$ and
$\mathcal{C}_\sigma^+\cap\left(\{0\}\times X\right)=\emptyset$.

\emph{Case 1}: $\mu_1^+(p)/f_\infty<\gamma<\mu_1^+(p)/f_0$.

In this case, we only need to show that
\begin{equation}
\left(\frac{\mu_1^+(p)}{\gamma f_\infty},\frac{\mu_1^+(p)}{\gamma
f_0}\right) \subseteq\left\{\mu\in\mathbb{R}\big|(\mu,u)\in \mathcal{C}_\sigma^+\right\}.\nonumber
\end{equation}
We divide the proof into two steps.

\emph{Step 1}: We show that if there exists a constant $M>0$ such that $\eta_n\subset(0,M]$ for $n\in \mathbb{N}$ large enough.

In this case it follows that $\Vert y_n\Vert\rightarrow+\infty.$
Similar to the proof of Theorem 4.1, we divide the equation
\begin{equation}
-\text{div}(\varphi_p(\nabla y_n))-\mu_n \gamma m(x)\varphi_p(y_n)=\mu_n \gamma m(x)\xi(y_n)\nonumber
\end{equation}
by $\Vert y_n\Vert_{C^{1,\alpha}(\overline{\Omega})}$ and set $\overline{y}_n =y_n/\Vert y_n\Vert_{C^{1,\alpha}(\overline{\Omega})}$. Since $\overline{y}_n$ is bounded in ${C^{1,\alpha}(\overline{\Omega})}$,
after taking a subsequence if
necessary, we have that $\overline{y}_n \rightharpoonup \overline{y}$ for some $\overline{y} \in {C^{1,\alpha}(\overline{\Omega})}$ and $\overline{y}_n
\rightarrow \overline{y}$ in ${C}(\overline{\Omega})$.
Using the similar method to the proof of Theorem 4.1, we can obtain
\begin{equation}
\lim_{n\rightarrow+\infty}\frac{ \xi(y_n(t))}{\Vert y_n\Vert_{C^{1,\alpha}(\overline{\Omega})}^{p-1}}=0\,\, \text{as\,\,} n\rightarrow+\infty.\nonumber
\end{equation}
By the compactness of $R_p:L^\infty(\Omega)\rightarrow X$ (see [\ref{DPM}]), we obtain
\begin{equation}
-\text{div}(\varphi_p(\nabla\overline{y}))-\left(\overline{\mu} m(x) \varphi_p(\overline{y})\right)=0,\nonumber
\end{equation}
where $\overline{\mu}=\underset{n\rightarrow+\infty}\lim\mu_n$, again choosing a subsequence and relabeling if necessary.
The rest proof of this step is the same as the proof of Theorem 4.1.

\emph{Step 2}: We show that there exists a constant $M$ such that $\mu_n
\in(0,M]$ for $n\in \mathbb{N}$ large enough.

On the contrary, we suppose that $\lim_{n\rightarrow +\infty}\mu_n=+\infty.$
Since  $\left(\mu_n, y_n\right) \in \mathcal{C}_\sigma^+$, it follows that
\begin{equation}
\text{div}\left(\varphi_p(\nabla y_n)\right)+\gamma\mu_n
m(x)\frac{f(y_n)}{\varphi(y_n)}\varphi(y_n)=0\,\, \text{in}\,\,\Omega^+,\nonumber
\end{equation}
where $\Omega^+=\{x\in\Omega\big|m(x)>0\}$. By Theorem 2.6 of [\ref{AY}], we have $y_n$ must change sign in $\Omega^+$, which contradicts Lemma 5.1.
The rest proof of is similar to the proof of Theorem 4.1.

\end{document}